\documentclass[30pt]{amsart}
\usepackage{epsfig, amssymb}
\newtheorem{theorem}{Theorem}[section]
\newtheorem{proposition}[theorem]{Proposition}
\newtheorem{lemma}[theorem]{Lemma}
\newtheorem{corollary}[theorem]{Corollary}

\theoremstyle{definition}

\theoremstyle{remark}
\newtheorem{remark}[theorem]{Remark}

\numberwithin{equation}{section}

\newcommand{\be}{\begin{equation}}
\newcommand{\ee}{\end{equation}}
\newcommand{\e}{\epsilon}

\newcommand{\vp}{\varphi}

\newcommand{\st}{\Xi}

\newcommand{\habr}{\hbar}
\newcommand{\eto}{T}

\begin{document}

\title{Semiclassical approximation and non-commutative geometry}
\vskip 1cm
\hfill \textit{\`a mon p\`ere}
\vskip 1cm
\author{ T. Paul\\ CNRS and CMLS, Ecole Polytechnique}
\maketitle
\vskip 0.3cm
\date{}
\vskip 1cm
\tableofcontents
\vskip 1cm

\begin{abstract}

We consider the long time semiclassical evolution for the linear Schr\"odinger equation. We show that, in the case of chaotic
underlying classical dynamics and for times up to $\hbar^{-2+\epsilon},\ \epsilon>0$, the symbol of a propagated observable by the
corresponding von Neumann-Heisenberg equation is, in a sense made precise below,  precisely obtained by the push-forward of the
symbol of the observable at time $t=0$. The corresponding definition of the symbol calls upon a kind of Toeplitz quantization 
framework, and the symbol itself is an element of  the noncommutative algebra of
the (strong) unstable foliation of the underlying dynamics.

\end{abstract}
\vskip 1cm
\section{Introduction}\label{intro}
In this note we  consider the long time  semiclassical evolution through the linear Schr\"odinger equation, or more precisely to the associated 
von Neumann equation
\be\label{von2}
i\hbar\frac d {dt}O^t=[O^t,H],
\ee
where $H$ is a Schr\"odinger operator $H=-\hbar^2\Delta+V$ with smooth confining  $V$  ($V(x)\to+\infty$ as
$\vert x\vert\to\infty$)- or a more general   semiclassical pseudodifferential operator of principal symbol $h$,
elliptic and selfadjoint on the Hilbert space $L^2(\mathcal M)$, where $\mathcal M$ is a manifold of dimension $n+1$.

It is well known \cite{BGP,BR} that, for times smaller than $C\log{\frac 1 \hbar}$, $C$ small enough, 
$O^t$ is still a Weyl (semiclassical) pseudodifferential operator and that its principal symbol is the push-forward 
of the initial one by the Hamiltonian flow associated to the principal symbol $h$ of $H$.
It is easy to get convinced  \cite{P1} that this is already not true for large values of $C$ (greater than $\frac 2 3$ 
times the natural Liapunov exponent of the flow).

Through this paper we will suppose that the Hamiltonian flow generated by $h$ is Anosov, 
and moreover that there exists a smooth action 
of $\mathbb R^{2(n+1)}$ on $T^*\mathcal M$, 
$(\mu,\nu;s,p)\in\mathbb R^{2(n+1)}\to
\st^\nu\circ\Lambda^\mu\circ T^{s,p}$, satisfying

\be\label{horf}
\st^\nu\circ\Lambda^\mu\circ T^{s,p}\circ\Phi^t=\Phi^t\circ\st^{e^{-\lambda t}\nu}\circ\Lambda^{e^{\lambda t}\mu}\circ T^{s+tp,p},\ \lambda>0.
\ee
(\ref{horf}) is obviously mimicked on the case of the geodesic flow on surface of constant curvature ($\Lambda^\mu$ (resp. $\st^\nu$) 
is the (resp. anti-)horocyclic flow,
 $T^{s,0}$ is the geodesic one and $T^{0,p}$ corresponds to a shift of energy) \cite{hor}, but we will not suppose that we are in this case and we'll 
use only (\ref{horf}). Moreover we will restrict this note to the bidimensional situation $n=1$ (the extension to any $n$,
keeping (\ref{horf}), is straightforward) and take $\mathcal M=\mathbb R^{n+1}$. The proofs are local and therefore are easily
adaptable to the non flat situation using the results of \cite{PU}. Finally we could extend some of our results to the case of variable Liapunov exponents.

 We will suppose that $O^{t=0}$ is a semiclassical pseudodifferential operator with smooth symbol supported in $h^{-1}(I)$ 
 for some interval $I\subset\mathbb R$ such that $h^{-1}(I)$ is compact.

 We first define, associated to $a\in\mathcal D(\mathbb R^{n+1}),\ \Vert a\Vert_{L^2}=1$, and the family of so-called (Gaussian) coherent states 
 $\vp_{(p,q)}(x):=(\pi\hbar)^{-\frac{n+1}4}e^{-\frac{(x-q)^2}{2\hbar}}e^{i\frac{px}\hbar}, \ (p,q)=z\in\mathbb R^{2(n+1)}$, 
 the family of Lagrangian states:
 \be\label{etat1}
\psi_z^a:=\int \exp{\big(\frac{i}\hbar\int\limits_z^{\Lambda^\mu\circ T^{s,0}(z)}\eta\big)}a(\mu,s)\vp_{\Lambda^\mu\circ T^{s,0}(z)}\frac{d\mu ds}{\hbar^{\frac{n+1}2}},\ \eta:=\xi.dx\
(\mbox{symplectic potential on}\ T^*\mathbb R^n).
\ee
It is easy to see that, microlocally in the interior of $I$ and for a support of $a$ small enough, the operator defined by 
$\int_{h^{-1}(I)}\vert\psi_z^a\rangle\langle\psi_z^a\vert \frac {d^{n+1}z}{\hbar^{n+1}}$ 
(here we denote by $\vert\psi\rangle\langle\psi\vert$ the operator:$\vp\to\langle\vp,\psi\rangle\psi$)
is equal to the identity modulo $\hbar^\infty$.

The key idea of this paper will be to write any pseudodifferential operator in the form
\be\label{pseud}
O=\int_{h^{-1}(I)}\vert\psi_z^{O_za}\rangle\langle\psi_z^a\vert \frac{d^{n+1}z}{\hbar^{n+1}}
\ee
for a suitable family of bounded pseudodifferential operators $O_z$.

The interest of such a formulation will be the fact that it is preserved by the evolution through (\ref{von2}). 
More precisely we prove in theorem \ref{1} that, for $n=1$ and
 any $0\leq t\leq C\hbar^{-2+\epsilon},\ \epsilon>0,$ there exists a bounded operator $O_z^t$ such that the solution $O(t)$ 
 of (\ref{von2}), $O(t=0)$ being
 microlocalized on $h^{-1}(I)$, satisfies

\be\label{oper}
\Vert \int_{h^{-1}(I)}\vert\psi_z^{O_z^ta}\rangle\langle\psi_z^a\vert \frac{d^{n+1}z}{\hbar^{n+1}}
-O(t)\Vert_{\mathcal B(L^2)}=O(\hbar^\infty)
\ee
(valid also for $O=$Identity$|_{h^{-1}(I)}$). This suggests to consider $O_z^t$ as the symbol of $O(t)$ at the point $z$.

In fact we will identify the symbol of $O(t)$ as a noncommutative object related to the space of leaves of the unstable 
foliation of the dynamics generated by the principal symbol $h$ of $H$.
Let us give the 
motivation behind this identification.

The classical limit of the equation (\ref{von2}) is the well known Liouville equation $\dot{\mathcal O}=\{\mathcal O,h\}$, 
where $\{.,.\}$ is the Poisson bracket on $T^*\mathcal M$, 
solved by the push forward of the initial condition by the Hamiltonian flow $\Phi^t$ associated to $h$. 
Though the flow is defined for all times $t$, the limit as $t\to-\infty$ doesn't have any meaning as a flow, 
being nevertheless the key ingredient of the theory of chaotic behaviour. In fact, if such a limit flow would exist it would be natural 
to say that it would be constant 
on the strong unstable manifold associated to any point $z\in T^*\mathcal M$ which, 
in our case, is the set of points $\Lambda_z=\{\Lambda^\mu\circ T^{s,0}(z),(\mu,s)\in\mathbb R^{n+1}\}$.
Therefore the pushforward of a smooth
 initial condition $\mathcal O$ by ``$\Phi^{-\infty}$" should be constant of each $\Lambda_z$, that is to say it should be a ``function" 
 on the space of leaves of the unstable foliation, orbits of the action of $\mathbb R^{n+1},\ (\mu,s)\to\Lambda^\mu\circ T^{s,0}$. 
 The leaves $\Lambda_z$ being usually dense on the energy shell, any (non constant) such function couldn't have 
 any regularity property (trace of the shearing off of the flow for long values of the time). The noncommutative geometry develops a topological theory for such singular spaces by, roughly speaking, replacing
 the algebra of continuous functions by a noncommutative one which, in the case of space of orbits of the action of a locally compact group, reduces to the crossed product
 of the algebra of continuous functions on the ambient manifold by the group. Let us note that this change of paradigm is invisible by the classical dynamics which is purely local.
 
 One of the main result of the present paper is to show that, approaching the limit $t\to -\infty$ on the time evolution of the classical dynamics 
 by a correlated semiclassical limit of the Schr\"odinger equation $\hbar\to 0$, $0<<t<\hbar^{-2+\epsilon}$, one recovers a dynamics 
 based on the noncommutative algebra of the strong unstable manifold, that is the ``space" of 
 the invariants of the local classical theory. 
 
 The noncommutative algebra of the unstable foliation is the geometrical setting of the classical limit
 of the long time quantum evolution. 
 
 Let us remark finally that 
  long time quantum 
 evolution creates  oscillations
  in the symbols of observables \cite{P1}.
 Therefore it is natural to consider the  microlocalization  of the symbol of the observable itself. 
 At the same time these oscillations are, 
 at each point of $T^*\mathcal M$, 
 along the unstable manifold, a
 highly non-linear object. It is then natural to expect that the good geometrical setting is not the cotangent bundle 
 over $T^*\mathcal M$ but precisely
 the unstable foliation, which is not a fibration in general, but which is an object  handleable by  noncommutative geometry.

 \vskip 1cm
\section{Propagation}\label{prop}
Let $O$ be a pseudodifferential operator whose symbol $\mathcal O$ is smooth and compactly supported inside $h^{-1}(I)$.
\vskip 1cm
\begin{theorem}\label{1}
Let us take $n=1$. There exists  bounded smooth and explicitly computable functions on $\mathbb R^{4}$, $\mathcal O^t_z\sim \mathcal O+\sum_{j=1}^\infty \mathcal O^t_j\hbar ^j$, 
such that, uniformly for 
$0\leq t\leq\hbar^{-2+\epsilon}$,
\[
\Vert e^{-i\frac{tH}\hbar}Oe^{+i\frac{tH}\hbar}-\int_{h^{-1}(I)}\vert\psi_z^{\widetilde {O_z}^ta}\rangle\langle\psi_z^a\vert
\frac{dz}{\hbar^2}\Vert_{\mathcal B(\mathcal H)}=O(\hbar^\infty),
\mbox{ $\widetilde {O_z}^t$ has total semiclassical Weyl symbol}
\]
\be\label{infty}
 \widetilde {\mathcal O^t_z}(\xi;x):=\mathcal O^t(\st^{e^{\lambda t}\nu}\circ\Lambda^{e^{-\lambda t}\mu}
\circ T^{s+tp,p}\circ\Phi^t(z)),\ x=(\mu,s),\ \xi=(\nu,p).
 \ee
 
\end{theorem}

 
\underline{Sketch of the proof:} the proof consist in several steps.
\begin{itemize}
\item we first prove the result for $t=0$. In order to do that we first show that, for $a$ with small enough support and microlocally on $h^{-1}(I')$ for
$I'\subset I$, $\int_{h^{-1}(I)}\vert\psi_z^{a}\rangle\langle\psi_z^a\vert
\frac{dz}{\hbar^n}=\mathbb I+O(\hbar^\infty)$, where $\mathbb I$ is the identity.
Since $\psi_z^a$ is a Lagrangian distribution, $O\psi_z^a=\psi_z^{a'}+O(\hbar^\infty)$ where $a'$ is obtained by the action of differentiable operators (transport equation).

\item $a$ having a frequency set included in the null section, one proves that $a'=\widetilde O^{0}_za+O(\hbar^\infty)$ 
where the Weyl symbol of $\widetilde O^{0}$
 has the form  $\mathcal O_z(\nu,p;\mu,s)\sim\sum\hbar^j\mathcal O_j(\st^{e^{\lambda t}\nu}\circ\Lambda^{e^{-\lambda t}\mu}
\circ T^{s,p}(z)))$.
 \item the next step is the heart of the proof. We want to show that $e^{-i\frac{tH}\hbar}Oe^{+i\frac{tH}\hbar}\psi_z^a=\psi_z^{a^t}$ 
 for some $a^t$ satisfying an equation that we derive, thanks to the main 
 hypothesis. In fact $a^t=\widetilde O^t_za$ where :
 \be\label{eq1}
 \dot{\widetilde O_z^t}=[\widetilde O_z^t,H_2+H_3] \mbox{, $H_2$ with a quadratic symbol and $H_3$ differential operator of third order}.
 \ee
 (\ref{eq1}) can be solved at any order,  $H^2$ being a quadratic operator (therefore giving an explicit solution) 
 and $H_3$ being treated by perturbation methods, after  microlocalizing near the zeroth section.
 \item thanks to this inoffensive microlocalization, we show that the preceding solution is valid 
 with an error term of the form
 \be\label{eq2}
 e^{-i\frac{tH}\hbar}Oe^{+i\frac{tH}\hbar}\psi_z^a=\psi_z^{a^t_k}+O((t\hbar^{2})^{k+1}\Vert a\Vert_{H^k(\mathbb R^n)}).
 \ee
 \item  taking $k>n$ and $t\leq\hbar^{-2+\epsilon}, \ \epsilon>0$ we get, since $h^{-1}(I)$ is compact, that
 $i\hbar\partial_tO^t_k=[O^t_k,H]+O_{\mathcal B(L^2)}(\hbar^{k\epsilon})$ where 
 $O^t_k:=\int_{h^{-1}(I)}\vert\psi_z^{\widetilde {a^t_k}}\rangle\langle\psi_z^a\vert
\frac{dz}{\hbar^n}$,  
from which we deduce, using the unitary of the propagator, 
$O^t=O^t_k+O(t\hbar^{k\epsilon-1})$ and, taking $k$ arbitrary, the result (\ref{infty}).
\end{itemize}
\begin{remark}\label{r}
As a Corollary of the proof of Theorem \ref{1} it is easy to prove that a similar result is still valid when we replace 
the functions $a$ by a 
$\hbar$-dependant ones of the form $a_\hbar(\dot)=\hbar^{-(n+1)\epsilon'/2}a(\hbar^{-\epsilon'}\dot)$ for $\epsilon'\geq 0$ small enough 
(see (\ref{eq2})).
More precisely, taking $\epsilon'<\epsilon$ (\ref{infty}) is still valid, and 
the conclusion of the last item is valid by replacing $\epsilon$ by $\epsilon-\epsilon'$. 
This allows to get the value $\e' =\frac 1 2$, that is the scaling corresponding to usual coherent states (take $a$ Gaussian), for times  up of the order $\hbar^{-\frac 3 2+\e},\  \e>0$.
\end{remark}
Let us mention  another Corollary of the proof of Theorem \ref{1}.
\vskip 1cm
\begin{proposition}\label{p}
There exist smooth functions $a^t\sim\sum_{j=0}^\infty a^t_j\hbar^j$ such that for $0\leq t\leq\hbar^{-2+\epsilon},\ \epsilon>0$,
\be\label{ppp1}
e^{+i\frac{tH}\hbar}\psi_z^a=e^{i\int_0^tpdq/\hbar}\psi_{\Phi^t(z)}^{\widetilde{a^t}}+O(\hbar^\infty) \mbox{ with }
\widetilde{a^t}(\mu,s)=e^{-it\hbar\partial^2_s/2}a^t(e^{\lambda t}\mu,s). 
\ee
In particular
\[
e^{+i\frac{tH}\hbar}\psi_z^a=e^{i\int_0^tpdq/\hbar}\psi_{\Phi^t(z)}^{\widehat{a^t}}+O(\hbar^\epsilon) \mbox{ with }
\widehat{a^t}(\mu,s)=e^{-it\hbar\partial^2_s/2}a(e^{\lambda t}\mu,s)\]
\end{proposition}
The preceding construction  is of course possible for $t\to -\infty$ verbatim by replacing the unstable 
by the stable foliation,  the flow $\Lambda^\mu$ by $\st^\mu$ and the states $\psi_z^a$ by $\vp_z^a$ constructed exactly the same way with $\st^\mu$ in place of $\Lambda^\mu$. 
Moreover the preceding estimations show that \eqref{ppp1} is also valid for $t=-\frac{\log{1/\hbar}}{2\lambda}:=t(\hbar)$. On can notice that 
$e^{-i\frac{t(\hbar)H}\hbar}\psi_z^a$ is, modulo $\hbar^\infty$ equal both to 
$e^{i\int_0^{-t(\hbar)}pdq/\hbar}\psi_{\Phi^{-t(\hbar)}(z)}^{\widetilde{a^{-t(\hbar)}}}$ and to
$e^{i\int_0^{t(\hbar)}pdq/\hbar}\vp_{\Phi^{t(\hbar)}(z)}^{\widetilde{b_\hbar}}$ since it is a semiclassical Hermite distribution associated to the isotropic manifold consisting in the trajectory issued from $\Phi^{t(\hbar)}(z)$. $\widetilde{b_\hbar}(\mu,s)=\hbar^{-\frac{n}4}b_\hbar(\mu/\sqrt\hbar,s)$ where $b_\hbar$ can be constructed as follows: $b_\hbar(.,s)=(M(s)a)(.,s)$ where $M(s)$ is the metaplectic quantization of the linear symplectic transformation mapping the tangent space of the unstable manifold at the point $\Phi^{-t(\hbar)+s}(z)$ to the tangent space of the stable manifold at the same point (\cite{P1}). Propagating $\vp_{\Phi^{t(\hbar)}(z)}^{\widetilde{b_\hbar}}$ by the equivalent of Proposition \ref{p} valid for negative times, we get (see also Remark \ref{r2} below) the
\begin{corollary}\label{corp}
There exist smooth functions $b^t\sim\widetilde{b_\hbar}+\sum_{j=1}^\infty b^t_j\hbar^j$ such that for times $t$ satisfying 
$-\hbar^{-\frac 3 2+\e}\leq t\leq-\frac{\log{1/\hbar}}{2\lambda},\ \epsilon>0$,
\be\label{corppp1}
e^{+i\frac{tH}\hbar}\psi_z^a=e^{i\int_0^tpdq/\hbar}\vp_{\Phi^t(z)}^{\widetilde{b^t}}+O(\hbar^\infty) \mbox{ with }
\widetilde{b^t}(\mu,s)=e^{-i(t-\frac{\log{1/\hbar}}{2\lambda})\hbar\partial^2_s/2}b^t(e^{\lambda t}\hbar^{1/2}\mu,s). 
\ee
In particular
\[
e^{+i\frac{tH}\hbar}\psi_z^a=e^{i\int_0^tpdq/\hbar}\vp_{\Phi^t(z)}^{\widehat{b^t}}+O(\hbar^\epsilon) \mbox{ with }
\widehat{b^t}(\mu,s)=e^{-i(t-\frac{\log{1/\hbar}}{2\lambda})\hbar\partial^2_s/2}\widetilde{b_\hbar}(e^{\lambda t}\hbar^{1/2}\mu,s).
\]
\end{corollary}

We see that, for $-\hbar^{-\frac 3 2+\e}\leq t\leq-\frac{\log{1/\hbar}}{2\lambda},\ \epsilon>0$, the propagator $e^{+i\frac{tH}\hbar}$ propagates a Lagrangian distribution associated to the unstable manifold of a given point $z$ to a Lagrangian distribution associated to the stable one of the same point.
This fact has to be compared to Fourier integral propagators which propagate Lagrangian manifold according to the underlying flow (or more general symplectic transformations): in the present case the stable manifold cannot be considered as the image of the unstable one by the  classical flow, but by a composition of flows with a ``symplectic mapping" responsible for the passage from the symbol $a$ to $\widetilde b_\hbar$ in \eqref{corppp1}. This fact can be compared to the case of the semiclassical propagation near an homoclinic trajectory studied in \cite{P3}. This induces a quantum mechanics inheritated interplay between the two classical Lagrangian invariant foliations. Let us mention finally that, for times as in Corollary \ref{corp}, the diagonal matrix elements of the propagator between  states $\psi_z^a$, key ingredients for the trace formula thanks to the decomposition of identity using the $\psi_z^a$s, will involve in the semiclassical limit the intersection of the stable and unstable manifolds of $z$, namely
 trajectories homoclinic to the one issued from $z$.
\begin{remark}\label{r2}
The same argument than  in Remark \ref{r} applies also to Proposition \ref{p}. This allows to compute the solution of equation \eqref{von2} up to times of the order $\hbar^{-\frac 3 2+\e},\  \e>0,$ in the framework of Toeplitz quantization. Namely writing the initial condition as :
\be\label{tq}
O=\int_{h^{-1}(I)}O(z)\vert\psi_z^{a}\rangle\langle\psi_z^a\vert \frac{d^{n+1}z}{\hbar^{n+1}}
\ee
with  $a(\eta)=(\pi\hbar)^{-(n+1)/4}e^{-\eta^2/2}$, we'll get for $0\leq t\leq \hbar^{-\frac 3 2-\e},\  \e>0$, thanks to Proposition \ref{p}, 

\be\label{tqp}
e^{-i\frac{tH}\hbar}Oe^{+i\frac{tH}\hbar}=\int_{h^{-1}(I)}O(z)\vert
\psi_{\Phi^t(z)}^{\widetilde{a^t}}\rangle\langle\psi_{\Phi^t(z)}^{\widetilde{a^t}}\vert
\frac{dz}{\hbar^{n+1}}+O(\hbar^\infty),
\ee 
with $\widetilde{a^t}$ given in Proposition \ref{p}. In other words we have the following
\begin{proposition}\label{tqpp}
\be\label{tqp2}
e^{-i\frac{tH}\hbar}Oe^{+i\frac{tH}\hbar}=\int_{h^{-1}(I)}O(\Phi^{-t}(z))\vert
\psi_{z}^{\widetilde{a^t}}\rangle\langle\psi_{z}^{\widetilde{a^t}}\vert
\frac{dz}{\hbar^{n+1}}+O(\hbar^\infty).
\ee 
\end{proposition}
Formula \eqref{tqp2} gives a kind of ``exact" quantum-classical evolution correspondence which might look close to the one expressed in 
Theorem \ref{1} since the ``symbol" $O(z)$ in \eqref{tqp2} is exactly propagated by the classical flow. Nevertheless Formula \eqref{tqp2} is highly nonlocal, in the sense that it will involve a nonlocal part of a function on which the right hand side of \eqref{tqp2} will apply, to be compared to  the right hand side of \eqref{infty}, except for times up to $(1-\e)\frac{\log{1/\hbar}}\lambda,\ \e>0,$ for which \eqref{tqp2} is a true alternative to Weyl quantization. Indeed using the results of \cite{BR} we know that for $0\leq t\leq (1-\e)\frac{\log{1/\hbar}}{2\lambda},\ \e>0,\  e^{-i\frac{tH}\hbar}Oe^{+i\frac{tH}\hbar}$ is an operator with principal Weyl symbol $O(\Phi^{-t}(z))$ which can be expressed in the form \eqref{tq}. It is easy to check that this leads to the
\begin{corollary}\label{cortqpp}
Let $0\leq t\leq (1-\e)\frac{\log{1/\hbar}}\lambda,\ \e>0.$ Then
\be\label{tqp3}
e^{-i\frac{tH}\hbar}Oe^{+i\frac{tH}\hbar}=\int_{h^{-1}(I)}O(\Phi^{-t}(z))\vert
\psi_{z}^{\widetilde{a}^{t-\frac{\log{1/\hbar}}{2\lambda}}}\rangle\langle\psi_{z}^{\widetilde{a}^{t-\frac{\log{1/\hbar}}{2\lambda}}}\vert
\frac{dz}{\hbar^{n+1}}+O(\hbar^\e).
\ee 
Note that $\psi_{z}^{\widetilde{a}^{t-\frac{\log{1/\hbar}}{2\lambda}}}$ are still microlocalized near $z$ and that subsymbols can also be computed.
\end{corollary}
Let us notice that Corollary \ref{cortqpp}, being purely local, can be proven in more general situations that the one in the present paper (see \cite{P1} for further details).

Another fact to be reported is that  the ``quantization process" in \eqref{tqp2}  and \eqref{tqp3}  depends on time, on the contrary to  the one in \eqref{infty} which doesn't, modulo the incorporation of the action of  the operator $\widetilde O_z^t$ on $a$ in $\psi_z^{\widetilde {O_z}^ta}$ which is reduced to $O(\Phi^{-t}(z))\times Identity$ in the framework of the (usual) Toeplitz quantization. 

Furthermore the next section will give a geometrical interpretation of Theorem \ref{1}.
\end{remark}

\section{Noncommutative geometry interpretation}\label{noncom}
We first prove the following Lemma.
\vskip 0.5cm
\begin{lemma}\label{deflemm}
Let us define
\be\label{symbole}
\sigma_{\mathcal O^t}(z,z'):=F(z_1,x,s;x'.s')
\ee where $z=\Lambda^xT^{s,0}(z_1),\ z'=\Lambda^{x'}T^{s',0}(z_1)$ and $F(z_1,.,.)$ is the integral symbol of an operator of Weyl symbol given by (\ref{infty}). Then
$\sigma_{\mathcal O^t}(z,z')$ doesn't depend on $z_1$.
\end{lemma}
The Lemma is easily proven by the translation invariance properties of the Weyl quantization procedure.

We want to identify $\sigma_{\mathcal O^t}$ as an element of the
crossed product $\mathcal A$ of the algebra $\mathcal C_I$ of continuous functions on $h^{-1}(I)$ by the group $\mathbb R^{n+1}$ under the action (\ref{horf}).
A function $\sigma(z,z'),\  z'\in\Lambda_z$ cam be seen as a continuous function from $G$ to $\mathcal A$ by 
\be
f(\mu,t)(z)=\sigma(z,z'),\ \ \ z'=\Lambda^\mu\circ T^t(z).
\ee
Moreover we get an action of $G$ on $\mathcal C_I$ by, $\forall g\in G$,
\be
\alpha_{(\mu,t)} h(z)=h(\Lambda^\mu\circ T^t(z)).
\ee
The algebra structure on $\mathcal C_I\rtimes_\alpha G$ is given by the $\star$-product $(f_1\star f_2)(g)=\int f_1(g_1)\alpha_{g_1}(f_2(g_1^{-1}g))dg_1$.

An easy computation, using Theorem \ref{1} and the symbolic property of Weyl quantization shows easily that, at leading order and for all $0\leq t_1,t_2\leq\hbar^{-2+\epsilon}$,
\be\label{symb1}
\sigma_{\mathcal O^{t_1} O^{t_2}}\sim \sigma_{\mathcal O^{t_1}}\star\sigma_{\mathcal O^{t_2}}.
\ee
Moreover the norm $\vert\vert\vert.\vert\vert\vert$ on $\mathcal A$ is equal to the supremum over $z$ of the operator norm on $L^2(\mathbb R^{n+1})$ of the operator 
of integral kernel $\sigma(z,z')$ (more precisely $\mathcal A\rtimes_\alpha G$ is the completion of the algebra of compactly supported kernels
$\sigma(z,z')$ with respect to the norm $\vert\vert\vert.\vert\vert\vert$.

We can also give a corresponding interpretation of the vectors $\psi_z^a$. 
Let us define $\alpha\in\mathcal A$ by, for 
$z'=\Lambda^\mu\circ T^{s,0}(z)$ , $\alpha(z,z'):=a(\mu,s)$. 
Then   $\psi_z^a=\psi^\alpha:=\int_{\Lambda_z} e^{\frac i\hbar\int_z^{z'}\eta}\alpha(z',z)\vp_{z'}dz'$.

We associate to any element $\gamma$ of  $\mathcal A$ an operator $\eto(\gamma)$ on $L^2(\mathbb R^{n+1})$ defined by
\be\label{def1}
\eto(\gamma):=\int_{h^{-1}(I)} \vert\psi^{\gamma\star\alpha}\rangle\langle\psi^\alpha\vert\frac{dz}{\hbar^{n+1}}.
\ee
In particular a bounded pseudodifferential operator is such an operator (with $\gamma\sim\sum_{j=0}^\infty\gamma_j\hbar^j$). 
Moreover, by definition of the norm $\vert\vert\vert.\vert\vert\vert$, $\eto(\gamma)$ is a bounded operator for all 
$\gamma\in \mathcal A\rtimes_\alpha G$ and it is easy to see, using arguments of the proof of theorem \ref{1}, 
that $\eto(\gamma)$ is bounded \textit{uniformly with} $\hbar\in [0,1]$ for $\gamma$ compactly supported. 
Noting that (\ref{def1}) is a way of writing (\ref{pseud}) we get:
\vskip 1cm
\begin{theorem}\label{2}
For $0\leq t\leq\hbar^{-2+\epsilon}$ there exist $\Gamma^t$ of ``symbol" $\gamma^t\sim\sum_{j=0}^\infty\gamma^t_j\habr^j\in\mathcal A$ such that
\be
\eto(\gamma)^t:=e^{-i\frac{tH}\hbar}\eto(\gamma) e^{+i\frac{tH}\hbar}=\eto(\gamma^t)+O(\hbar^\infty)\ \ \ \ \ \mbox{ with }
\ \ \ \ \gamma^t_0=\Phi^t\#\gamma_0+O(\hbar^\epsilon).
\ee
Moreover the leading order symbol of $\eto(\gamma)^{t_1}\eto(\gamma)^{t_2}$ is $\gamma^{t_1}\star\gamma^{t_2}$.
\end{theorem}
\underline{Sketch of the proof:} 
\begin{itemize}
\item the fact that Theorem \ref{1} is valid also for operators defined by (\ref{def1}) is contained in the proof of Theorem \ref{1} itself.
\item the fact that the symbol of $\eto(\gamma^t)$ is in the completion by the norm $\vert\vert\vert.\vert\vert\vert$ is obtained by the Calderon-Vaillancourt theorem, since
$\mathcal O(\st^{e^{\lambda t}\xi}\circ\Lambda^{e^{-\lambda t}x}
\circ T^{s+tp+t,p}(z))$ is bounded and smooth, therefore defines a bounded (non semiclassical) pseudodifferential operator.
\item the product formula of principal symbols is nothing but (\ref{symb1}).

\end{itemize}
Let us remark also that an extension on the lines of Remark \ref{r} is also valid in this framework.
\vskip 1cm
\section{Semiclassical measures}\label{semi}

In the same way that one associates to a vector $\psi$ (or density matrix) 
the quantity $\langle\psi,\vp_z\rangle\vert^2$ considered as a measure  by the formula
$\langle\psi,O\psi\rangle=\int \mathcal O_T(z,\bar z) \vert\langle\psi,\vp_z\rangle\vert^2dz$, 
where $\mathcal O_T$ is the Toeplitz symbol of $O$, on can associate to $\psi$ (or a density matrix) a sort 
of ``off-diagonal" version by the quantity $R_\psi(z,\bar{z'}):=\langle \vp_z,\psi\rangle\langle\psi,\vp_{z'}\rangle$ for $z'\in\Lambda_z$. 

$R_\psi$ can be considered as an element of the dual of a (dense) subalgebra of $\mathcal A$ and will have better properties of semiclassical 
propagation. For sake of shortness we express the result in the case of eigenvectors of the Hamiltonian $H$, leaving the 
straightforward derivation
 for $R_{e^i\frac{tH}\hbar}\psi$ in the same topology.
\vskip 1cm
\begin{theorem}\label{3}
Let us define for $I$ compact interval of $\mathbb R$, $\mathcal D_I\sim\mathcal D(h^{-1}(I)\times \mathbb R^n)$ the subalgebra of 
smooth compactly supported elements of $\mathcal A$. 
Let $\psi$ be an eigenfunction of $H$. 

Then, restricted to $z\in h^{-1}(I)$, $R_\psi(z,\bar{z'})$ considered as a function 
on $h^{-1}(I)\times \mathbb R^n$, belongs to $\mathcal D_I'$. Moreover, in the weak-* topology and, $\forall \epsilon >0$, 
uniformly for $0\leq t\leq \hbar ^{-2+\epsilon}$
\be\label{pf}
\Phi^t\# R(z,\bar{z'})
=R(z,\bar{z'})+O(\hbar^\epsilon).
\ee

\end{theorem}
The proof consists in using $\Gamma$, a quantization using $a_\hbar$ as in Remark \ref{r} of a symbol $\gamma(z,z')$ 
compactly supported in
$z'$. Writing the formula for $\langle \vp_{z},\Gamma \vp_{z'}\rangle$ we get that $R_\psi\in\mathcal D_I'$, and it is easy to see that
Theorem \ref{1} applies to $\Gamma$, out of which we derive (\ref{pf}).

\vskip 1cm
\section{Perspectives}\label{persp}
Other situations with a noncommutative semiclassical limit can be treated, e.g. the integrable cases.

In this paper we presented only preliminary results concerning the quantization of algebra of the unstable foliation. In particular more symbolic results can be obtained in full generality.

The construction of Section \ref{semi} is of course possible for $t\to -\infty$ verbatim by replacing the unstable 
by the stable foliation, and the flow $\Lambda^\mu$ by $\st^\nu$. We believe that it is
 possible to construct operators whose (noncommutative) symbols will be concentrated on the intersection of the two foliation, and to derive a  result similar to
 Theorem \ref{3} by some invariance property along homoclinic trajectories. All these works are in progress.

\end{document}